\documentclass{article}

\usepackage{graphicx}

\usepackage{amscd}

\usepackage{amsmath}

\usepackage{amsfonts}

\usepackage{amssymb}

\begin{document}

\centerline {\large {\bf C*-algebras of commuting endomorphisms}}

\bigskip

\centerline{\tt by VALENTIN DEACONU}

\bigskip

\hspace{70mm}{\em dedicated to \c{S}erban Str\u{a}til\u{a} on his
$60^{th}$ birthday}

\bigskip

\begin{abstract}Given a compact space $X$ and two commuting continuous open surjective maps
$\sigma_1,\sigma_2:X\rightarrow X$, we construct certain
C*-algebras that reflect the dynamics of the ${\mathbb
N}^2$-action. When the maps $\sigma_1, \sigma_2$ are local
homeomorphisms, these are groupoid algebras, but in general, we
will use a Cuntz-Pimsner algebra associated to a product system of
Hilbert bimodules in the sense of Fowler.

The motivating example for our construction is the dynamical system associated with
a rank two graph by Kumjian and Pask. We consider also
a two-dimensional subshift of Ledrappier, the case of two covering maps of the circle,
and the two-dimensional Bernoulli shift.\end{abstract}

\bigskip

{\bf \S 1. Introduction}

\bigskip

There are different approaches to define the crossed product of a
C*-algebra by a semigroup of endomorphisms in the works of G.
Murphy, I. Raeburn, M. Laca and others (see \cite{M1}, \cite{M2},
\cite{LR}, \cite{LiR}, \cite{F}). They are using covariant
isometric or partial isometric representations of the positive
cone in a totally ordered abelian group,  are realizing the
crossed product by a cancellative abelian semigroup of injective
endomorphisms as a corner of the covariance algebra of a classical
C*-dynamical system, or are using a family of Hilbert bimodules
indexed by the quasi-lattice ordered semigroup and generalized
Cuntz-Pimsner algebras .

Our examples are mainly related to an abelian ${\mathbb N}^2$-
dynamical system, and we will consider in section 2 a groupoid
approach, similar to \cite{D}, and  a Hilbert bimodule approach in
the next section. We will see how the definition of the inner
product, right and left actions will affect the outcome, and how
the dynamics of the ${\mathbb N}^2$-action is reflected in the
corresponding C*-algebra. We define the inner products by using
the transfer operators associated to some conditional expectations
onto the ranges of the endomorphisms as in \cite{E1},\cite{E2}. In
the case of local homeomorphisms, the resulting Hilbert bimodules
are finitely generated and projective, therefore the left action
consists entirely of compact operators.  In general, the transfer
operators are defined using a family of measures as in \cite{Br},
and it may happen that the left action has trivial intersection
with the compacts. For two commuting endomorphisms, we will
construct a product system of Hilbert bimodules over the semigroup
${\mathbb N}^2$(see \cite{F}), and the corresponding Cuntz-Pimsner
algebra will play the role of the groupoid algebra. Sometimes,
this algebra could be understood as an iterated Pimsner
construction, by extending the scalars and by using the universal
properties.
\bigskip

{\bf \S 2. The groupoid approach}

\bigskip

Consider a compact  space $X$ with two commuting continuous open
surjections
 $\sigma_1,\sigma_2:X
\rightarrow X.$ In analogy to the groupoids associated to a
covering map (see \cite{D}), we define the equivalence relations
\[R_n=\{(x,y)\in X\times X\;\mid\; \sigma^nx=\sigma^ny\}
\;\;{\mbox{and}\;\;} R=\bigcup_{n\in {\mathbb N}^2} R_n.\] Here
$n=(n_1,n_2)$ and $\sigma^n=\sigma_1^{n_1}\sigma_2^{n_2}$. We
consider ${\mathbb N}^2$ directed by the partial order
\[(n_1,n_2)\leq (m_1,m_2)\; {\mbox{if}}\; n_i\leq m_i, i=1,2.\]

We put the induced topology on $R_n$, and the inductive limit
topology on $R$.

\bigskip

{\bf Proposition 1}. If $\sigma_i$ are local homeomorphisms, then
all $R_n$ are r-discrete,  $C^*(R)$ is well defined, and
\[C^*(R)=\lim_{\rightarrow}C^*(R_n).\]

\bigskip

{\em Proof}.  See pages 122-123 in \cite{Re}.

\bigskip

 Let $\Gamma=\Gamma(\sigma_1,\sigma_2)=\{(x,p-q,y)\in
X\times {\mathbb Z}^2\times X\;\mid\; \sigma^px=\sigma^qy\}.$ With
the usual operations,
\[(x,k,y)\cdot(y,l,z)=(x,k+l,z),\;\; (x,k,y)^{-1}=(y,-k,x),\]$\Gamma$ is a
groupoid with the unit space identified with $X$. The isotropy
group bundle is
\[I=\{(x,p-q,x)\in X\times {\mathbb Z}^2 \times X\mid \;\sigma ^px=\sigma ^q x\}.\]

\bigskip

{\bf Proposition 2}. If all the equivalence relations $R_n$ are
r-discrete, then $\Gamma$ is a locally compact Hausdorff groupoid
with a Haar system.

\bigskip

{\em Proof}. The topology of $\Gamma$ is defined by the cylinder
sets
\[Z(U,V,p,q)=\{(x,p-q,(\sigma\mid_V)^{-q}\circ \sigma^p(x)), x\in U\},\]
where $p,q\in{\mathbb N}^2$, and $U$ and $V$ are open sets of $X$
such that $\sigma^p\mid_U$ and $\sigma^q\mid_V$ are
homeomorphisms. A Haar system on $R$ could be extended to
$\Gamma$.

\bigskip

{\bf Remark 1}. Assuming that $C^*(\Gamma)$ is well defined (which
certainly happens in the case $\sigma_i$ are local
homeomorphisms), the cocycle
\[c:\Gamma\rightarrow {\bf Z}^2, c(x,n,y)=n,\]
 induces a ${\mathbb T}^2$-action on $C^*(\Gamma)$, via
\[(z\cdot f)(x,n,y)=z^n\cdot f(x,n,y), f\in C_c(\Gamma), z\in {\mathbb
T}^2.\]

\bigskip

{\bf Remark 2}. If $\sigma_i$ are homeomorphisms, then
  $C^*(\Gamma)$ is isomorphic
to the usual crossed product $C(X)\times {\mathbb Z}^2$. In
general, we will see that $C^*(\Gamma)$ is a certain crossed
product of $C(X)$ by the semigroup ${\mathbb N}^2$.

\bigskip

{\bf Definition}. The {\em orbit} of $x\in X$ is defined to be
\[ O(x)= \bigcup _{k\in{\mathbb N}^2} \sigma ^{-k}(\sigma ^kx), \] where
$\sigma ^{-k}y= \{z\in X | \sigma ^kz= y\}.$

\noindent We say that $\sigma$ is {\em minimal} if each orbit is
dense and that $\sigma$  is {\em essentially free} if \[\{x\in X
\mid \, \forall k,l ,\,\sigma ^kx = \sigma^lx\,\Rightarrow \,k=l
\}\] is dense in $X$.

\bigskip

{\bf Proposition 3}. If $\sigma$ is minimal, then $C^*(R)$ is
simple. Moreover, if $\sigma$ is essentially free, then $
C^*(\Gamma )$ is also simple.

\bigskip

{\em Proof.} Since $\sigma$ is minimal and
\[ O(x) = \{y \in X \mid (x,y)\in R \},\]
it follows that there are no nontrivial open invariant subsets. If
$\sigma$ is also essentially free, then the groupoid $\Gamma$ is
essentially principal in the sense of Renault (see definition
II.4.3 of \cite {Re}). We can  apply proposition II.4.6 of \cite
{Re}, where the ideals of an essentially principal groupoid are
characterized.

\bigskip

{\bf Example 1}. Consider $X$ the infinite path space of a rank
two graph in the sense of Kumjian and Pask, and $\sigma_1,
\sigma_2$ the horizontal and vertical shifts, which are local
homeomorphisms (see \cite {KP}). For simplicity, we describe $X$
in a particular situation. Start with two finite graphs $(G_1,V)$
and $(G_2,V)$ with the same set of vertices $V$, such that there
is a bijection $\rho: G_1\ast G_2\rightarrow G_2\ast G_1$. Here
\[G_1\ast G_2=\{(e,f)\in G_1\times G_2\mid\; s(e)=r(f)\},\]
where $s$ and $r$ are the source and range maps. For example, this
happens if the vertex matrices commute. Then $X$ could be thought
as an infinite grid in the first quadrant, where the vertices at
the lattice points are joined by horizontal edges (from $G_1$) and
vertical edges (from $G_2$). The unique factorization property
translates in the fact that the bijection $\rho$ determines the
other two edges of each unit square, once we fix a horizontal
edge, followed by a vertical edge (see section 6 in \cite{KP}). It
follows that any finite rectangular grid is determined by one
horizontal side $e_1e_2...e_m$ and one vertical side
$f_1f_2...f_n$. This ensures that  $\sigma_i$ are local
homeomorphisms of the Cantor set $X$.

The C*-algebra of the corresponding groupoid $\Gamma$ is strongly
Morita equivalent to a crossed product of an AF-algebra by the
group ${\mathbb Z}^2$. Under some mild conditions, the groupoid is
essentially free, and the C*-algebra is simple and purely infinite
(see section 4 in \cite{KP}).

\bigskip

{\bf Example 2}. Let $X$ be the unit circle ${\mathbb T}$, and
$\sigma_1, \sigma_2$ two covering maps, of indices $p_1$ and
$p_2$, respectively, with $|p_i|\geq 2$. In this case, each
$C^*(R_n)$ is of the form $C({\mathbb T})\otimes M_k$ for some
positive integer $k=k(n)$,  $C^*(R)$ is a circle algebra, and
$C^*(\Gamma)$ is simple and purely infinite, since the orbits are
dense, and the groupoid is essentially free and locally
contracting (see \cite{A-D}).

\bigskip

{\bf Example 3} (Ledrappier). Let  $X\subset \{0,1\}^{{\mathbb
N}^2}$ be the subgroup defined by
\[x_{i+1,j}+x_{i,j}+x_{i,j+1}=0\;\;{\mbox{mod 2}}.\] The admissible patterns are
\[P=\left\{\begin{array}{cc}0&{}\\0&0\end{array},\;\;\begin{array}{cc}1&{}\\0&1\end{array},\;\;
\begin{array}{cc}1&{}\\1&0\end{array},\;\;\begin{array}{cc}0&{}\\1&1\end{array}\right\}.\]
The maps $\sigma_1$ and $\sigma_2$ are the horizontal and the
vertical shifts, which are local homeomorphisms in this case.
Indeed, let $x\in X$ be fixed, and consider the open neighborhood
of $x$, defined by
\[U=\{y\in X\;\mid\; y(0,0)=x(0,0)\}.\] Then $\sigma_1(U)=X$, and if $y,z\in U$ with $\sigma_1(y)=\sigma_1(z)$,
then $y(i,j)=z(i,j)$ for all $i\geq 1, \;j\geq 0$. But also
$y(0,0)=z(0,0)=x(0,0)$, and  since the only admissible patterns
are the ones in $P$, we get $y(0,1)=z(0,1)$. By induction, we get
$y(0,j)=z(0,j)$ for all $j\geq 0$, therefore $y=z$ and
$\sigma_1\mid_{U}$ is one to one. Similarly,  $\sigma_2$ is a
local homeomorphism.

The C*-algebra $C^*(R)$ is isomorphic to the CAR algebra
$UHF(2^{\infty})$, and $C^*(\Gamma)$ is again strongly Morita
equivalent to a crossed product of an AF-algebra by ${\mathbb
Z}^2$, like in the case of a rank 2 graph.

More general, we may consider  $V$ to be a finite alphabet and a
closed subset $X\subset V^{{\mathbb N}^2}$ (in the product
topology) that is $\sigma$-invariant, where \[\sigma^l
(x)(k)=x(k+l), \;x\in X, \;k,l\in {\mathbb N}^2.\] Such a
dynamical system $(X,\sigma)$ is a (two-dimensional) {\em Markov
shift} or a {\em subshift of finite type} if there exists a finite
set $F\subset {\mathbb N}^2$ and a set of {\em admissible
patterns} $P\subset V^F$ such that
\[X=X(F,P)=\{x\in V^{{\mathbb N}^2} \; \mid \;
\pi_F(\sigma^mx)\in P\; \mbox{for every }\; m\in {\mathbb
N}^2\},\] where $\pi_F$ is the projection onto $V^F$. If the shift
maps are local homeomorphisms, then the C*-algebra $C^*(\Gamma)$
will be strongly Morita equivalent to a crossed product of an
AF-algebra by ${\mathbb Z}^2$. But, in general, the shift maps are
not local homeomorphisms, as it can be seen in the case of the
full shift $X=\{0,1\}^{{\mathbb N}^2}$. The equivalence classes in
the equivalence relations $R_n$ are not discrete anymore, in fact
they are Cantor sets. For other interesting examples of higher
dimensional subshifts of finite type, see \cite{S}.

\bigskip

{\bf \S 3. The Hilbert bimodule approach}

\bigskip

For a C*-algebra $A$ and an injective unital endomorphism
$\alpha\in End(A)$ such that there is a conditional expectation
${\cal P}$ onto the range $\alpha(A)$, one can define a Hilbert
bimodule $E=A(\alpha, {\cal P})$, using the transfer operator
${\cal L}=\alpha^{-1}\circ {\cal P}$ (see \cite{E1}, \cite{E2}).
We start with the vector space $A$, and define the inner product,
right and left multiplications by the equations
\[<\xi, \eta>={\cal L}(\xi^*\eta), \;\;
\xi\cdot a=\xi\alpha(a),\;\;a\cdot \xi=a\xi.\] Then we have
\[<\xi, \eta\cdot a>=<\xi,\eta\alpha(a)>={\cal
L}(\xi^*\eta\alpha(a))=\alpha^{-1}({\cal
P}(\xi^*\eta\alpha(a)))=\]\[=\alpha^{-1}({\cal
P}(\xi^*\eta)\alpha(a))= \alpha^{-1}({\cal
P}(\xi^*\eta))a=<\xi,\eta>a.\] Using the Pimsner construction (see
\cite{P}), we get a Cuntz-Pimsner algebra ${\cal O}_{A(\alpha,
{\cal P})}$. Recall

\bigskip

{\bf Definition}. A {\em Toeplitz representation} of a Hilbert
$A$-bimodule $E$ in a C*-algebra $B$ is a pair $(\psi, \pi)$ with
$\psi :E\rightarrow B$ a linear map and $\pi :A\rightarrow B$ a
homomorphism, such that
\[\psi(\xi\cdot a)=\psi(\xi)\pi(a)\]
\[\psi(\xi)^*\psi(\eta)=\pi(<\xi, \eta>)\]
\[\psi(a\cdot\xi)=\pi(a)\psi(\xi).\]
The corresponding universal C*-algebra is called the Toeplitz
algebra of $E$, denoted by ${\cal T}_E$.
 There is a homomorphism
$\pi^{(1)}:{\cal K}(E)\rightarrow B$ which satisfies
\[\pi^{(1)}(\Theta_{\xi,\eta})=\psi(\xi)\psi(\eta)^*.\]
We say that $(\psi, \pi)$ is {\em Cuntz-Pimsner covariant} if
\[\pi^{(1)}(\phi(a))=\pi(a)\; \forall\;\; a\in\phi^{-1}({\cal K}(E)),\]
where $\phi:A\rightarrow L(E)$ is the left action. The
Cuntz-Pimsner algebra ${\cal O}_E$ is universal for Toeplitz
representations which are Cuntz-Pimsner covariant. There is a
gauge action of the circle group ${\mathbb T}$ on ${\cal O}_E$,
and the fixed point algebra  is denoted by ${\cal F}_E$.

\bigskip

{\bf Example 1}. For $A=C(X)$ and $\alpha$ induced by a local
homeomorphism $\sigma :X\rightarrow X$, we can take
\[({\cal P}f)(x)=\frac{1}{\nu(x)}\sum_{\sigma(y)=\sigma(x)}f(y),\]
where $\nu(x)$ is the number of elements in the fiber
$\sigma^{-1}(x)$. The corresponding algebra ${\cal
O}_{A(\alpha,{\cal P})}$ is isomorphic to $C^*(\Gamma(\sigma))$,
where $\Gamma(\sigma)$ is the Renault groupoid as in \cite{D},
\[\Gamma(\sigma)=\{(x,p-q,y)\in X\times {\mathbb
Z}\times X\;\mid\; \sigma^px=\sigma^qy\}.\] Indeed, it is known
that ${\cal K}(E)=E\otimes_A E^*$,  where $E^*$ is the adjoint of
$E$. That means $E^*=C(X)$ with the conjugate multiplication of
scalars, and with the left and right actions of $A$ interchanged.
Since $\xi\cdot f\otimes\eta^* =\xi\otimes \phi(f)\eta^*$, it
follows that
$\xi(x)f(\sigma(x))\eta^*(y)=\xi(x)f(\sigma(y))\eta^*(y)$ for all
$(x,y)\in X\times X$ and $f\in C(X)$. Hence $\sigma(x)=\sigma(y)$,
and  ${\cal K}(E)= {C(R_1)}$ as sets. We have
\[\Theta_{\xi,\eta}\Theta_{\xi',\eta'}(\zeta)=\Theta_{\xi,\eta}(\xi'<\eta',\zeta>)=
\xi<\eta,\xi'<\eta',\zeta>>=\xi<\eta,\xi'><\eta',\zeta>=\Theta_{\xi<\eta,\xi'>,\eta'}(\zeta),\]
therefore
\[(\xi\otimes\eta^*)(\xi'\otimes\eta'^*)(x,y)=\xi(x)<\eta,\xi'>(\sigma(x))\eta'(y)=
\sum_{\sigma(z)=\sigma(x)}\xi(x)\eta^*(z)\xi'(z)\eta'^*(y),\]and
the multiplication of compact operators is exactly the convolution
product on $C(R_1)$. Hence, ${\cal K}(E)=C^*(R_1)$ as
$C^*$-algebras. In the same way, using the fact that ${\cal
K}(E^{\otimes n})=(E^{\otimes n})\otimes_A(E^{\otimes n})^*$, we
get ${\cal K}(E^{\otimes n})=C^*(R_n)$. Taking inductive limits,
it follows that ${\cal F}_E=C^*(R)$, since in this case ${\cal
F}_E$ is generated by all ${\cal K}(E^{\otimes n}),\; n\geq 0$.
The isometry $v$ which induces the corner endomorphism of $C^*(R)$
is induced by the function $\gamma(x)=1$, regarded as an element
of $E$. Indeed, $<\gamma,\gamma>=1$, therefore $\psi(\gamma)$ is
an isometry in any Cuntz-Pimsner representation. The endomorphism
of ${\cal F}_E$ is given by the formula
\[(\xi_1\otimes\xi_2\otimes...\otimes\xi_n)
\otimes(\eta_1\otimes\eta_2\otimes...\otimes\eta_n)^*\mapsto
(\gamma\otimes\xi_1\otimes...\otimes\xi_n)\otimes(\gamma\otimes\eta_1\otimes...\otimes\eta_n)^*.\]

\bigskip

{\bf Remark}. The above Cuntz-Pimsner algebra ${\cal
O}_{A(\alpha,{\cal P})}$ could also be described as the universal
C*-algebra generated by a copy of $C(X)$ and an isometry $S$
subject to the relations
\[(i)\;\;Sf=\alpha(f)S,\;\;(ii) \;\;S^*fS={\cal L}(f),\;\; (iii) \;\; 1=\sum_{i=1}^mu_iSS^*u_i^*,\]
for all $f\in C(X)$ where
\[{\cal L}f(x)=\frac{1}{\nu(x)}\sum_{\sigma(y)=x}f(y),\]
and $\{u_1,u_2,...,u_m\}\subset C(X)$ such that $\displaystyle
f=\sum_{i=1}^mu_i{\cal P}(u_i^*f)$ for all $f\in C(X)$ (Theorem
9.2. in \cite{EV}). The corresponding Toeplitz algebra ${\cal
T}_{A(\alpha,{\cal P})}$ satisfies only the first two relations.

{\bf Note}. If $\sigma:X\rightarrow X$ is continuous, open and
surjective, but not a local homeomorphism, the  groupoid
construction  from \cite{D} may fail. The equivalence relation $R$
does not have an obvious Haar system, since the Haar systems on
$R_n$ may not be compatible. Since $R_n$ is not necessarily an
open subset of $R_m$ for $n\leq m$,  a continuous function on
$R_n$ has no natural extension to a function on $R_m$.

{\bf Example 2}. Let $\displaystyle X=\prod_1^\infty [0,1]$ and
$\sigma(x_1,x_2,...)=(x_2,x_3,...)$. Then $\sigma$ is not a local
homeomorphism, in fact the fibers are homeomorphic to the interval
$[0,1].$ Also, $R_n\subset R_{n+1}$ is not open, therefore we can
not conclude that a  continuous function on $R_n$ is naturally
extended to a continuous function on $R_{n+1}$.  We try to
overcome the lack of an obvious groupoid by using the existence of
certain families of measures on the fibers of $\sigma$ which will
define a transfer operator and a Hilbert bimodule over $C(X)$.

\bigskip

{\bf Definition}. Suppose that $\pi:X\rightarrow Y$ is a
continuous open surjection between locally compact Hausdorff
spaces. A family $\lambda =\{\lambda_y\}_{y\in Y}$ of positive
Radon measures on $X$ is called a $\pi$-system if the support of
$\lambda_y$ is contained in $\pi^{-1}(y)$ for each $y\in Y$, and
for each $f\in C_c(X)$, the function
\[\lambda(f)(y):=\int f(x)d\lambda_y(x)\]
lies in $C_c(Y)$. We say that a $\pi$-system is full if the
support of each $\lambda_y$ is all of $\pi^{-1}(y)$. The existence
of a full $\pi$-system is proved in (\cite{Bl}, Th\'eor\`eme 3.3)
for $X$ separable and $Y$ second countable.

For the above example, one can take the usual Lebesgue measure on
$[0,1]$ for each fiber. We get a conditional expectation, a
transfer operator, and a Hilbert bimodule $E$. The corresponding
Cuntz-Pimsner algebra is isomorphic to the Toeplitz algebra ${\cal
T}_{E}$ since $\phi^{-1}({\cal K}(E))=0$ in this case (see 3.3 in
\cite{Sc} and \cite{K}).

It was observed by Raeburn and Sims (\cite{RaS}) that the
C*-algebras associated to higher rank graphs could be obtained
from a product system of Hilbert bimodules in the sense of Fowler
(\cite{F}). We will do something similar for two commuting
endomorphisms of a C*-algebra $A$. Recall

\bigskip

{\bf Definition}(Fowler). Given $Q$ a countable semigroup with
identity $e$ and $p:E\rightarrow Q$ a family of Hilbert bimodules
over $A$, we say that $E$ is a {\em product system} over $Q$ if
$E$ is a semigroup, $p$ is a morphism, and for each $s,t\in
Q\setminus\{e\}$, the map $(\xi,\eta)\in E_s\times E_t\mapsto
\xi\eta\in E_{st}$ extends to an isomorphism between
$E_s\otimes_AE_t$ and $E_{st}$. We require that $E_e=A$ (with
$\phi_e(a)b=ab$), that the multiplications $E_e\times
E_s\rightarrow E_s$ and $E_s\times E_e\rightarrow E_s$ satisfy
$a\xi=\phi_s(a)\xi,\;\; \xi a=\xi\cdot a$, and that each $E_s$ is
essential as left $A$-modules, in the sense that $E_s$ is the
closed span of the elements $\phi_s(a)\xi$ with $a\in A$, $\xi\in
E_s$, and $\phi_s:A\rightarrow L(E_s)$ defining the left
multiplication.

\bigskip

{\bf Definition}. A {\em Toeplitz representation} of the product
system $E$ in a C*-algebra $B$ is a map $\psi:E\rightarrow B$ such
that

(1) For each $s\in Q, (\psi_s,\psi_e)$ is a Toeplitz
representation of $E_s$, where $\psi_s=\psi\mid_{E_s}$, and

(2) $\psi(\xi\eta)=\psi(\xi)\psi(\eta)$ for $\xi,\eta\in E$.

\noindent If in addition each $(\psi_s,\psi_e)$ is Cuntz-Pimsner
covariant, then $\psi$ is a {\em Cuntz-Pimsner representation}.

\bigskip

It was proved by Fowler (Propositions 2.8, 2.9 \cite{F}) that the
Toeplitz algebra ${\cal T}_E$ and  Cuntz-Pimsner algebra ${\cal
O}_E$ exist and are unique up to isomorphism.

Given $\alpha_i, i=1,2$ two commuting, injective unital
endomorphisms of a C*-algebra $A$ such that there exist commuting
conditional expectations ${\cal P}_i$ onto the ranges
$\alpha_i(A)$, we will construct a product system of Hilbert
bimodules over the semigroup ${\mathbb N}^2$. We take
$E_{(0,0)}=A$, $E_{(1,0)}=A(\alpha_1,{\cal P}_1),
E_{(0,1)}=A(\alpha_2,{\cal P}_2)$, and
$E_{(m,n)}=E_{(1,0)}^{\otimes m}\otimes_A E_{(0,1)}^{\otimes n}$
for all $(m,n)\in {\mathbb N}^2$. The semigroup structure is given
by the tensor product.

\bigskip

{\bf Lemma}. We have $A(\alpha_1,{\cal P}_1)\otimes_A
A(\alpha_2,{\cal P}_2)\simeq A(\alpha_1\circ\alpha_2,{\cal
P}_1\circ{\cal P}_2)$.

\bigskip

{\em Proof}. The map $\Phi:A(\alpha_1,{\cal P}_1)\otimes_A
A(\alpha_2,{\cal P}_2)\rightarrow  A(\alpha_1\circ\alpha_2,{\cal
P}_1\circ{\cal P}_2),
\Phi(\xi_1\otimes\xi_2)=\xi_1\alpha_1(\xi_2)$ induces the Hilbert
bimodules isomorphism. The inverse is $\Phi^{-1}(\xi)=\xi\otimes
1$.

\bigskip

{\bf Proposition}. If $\sigma_i:X\rightarrow X, i=1,2$ are two
commuting local homeomorphisms of a compact space $X$, and ${\cal
P}_i$ are given by
\[({\cal P}_if)(x)=\frac{1}{\nu_i(x)}\sum_{\sigma_i(y)=\sigma_i(x)}f(y),\]
where $\nu_i(x)$ is the number of elements in the fiber
$\sigma_i^{-1}(x)$, then ${\cal P}_1\circ{\cal P}_2={\cal
P}_2\circ{\cal P}_1$, and the Cuntz-Pimsner algebra ${\cal O}_E$
associated to the above product system is isomorphic to the
groupoid algebra $C^*(\Gamma(\sigma_1,\sigma_2))$ considered in \S
2.

\bigskip

{\bf Remark}. Such a product system over ${\mathbb N}^2$ could be
constructed from two more general Hilbert bimodules $E_1,E_2$ over
a C*-algebra $A$, such that $E_1\otimes_AE_2\simeq
E_2\otimes_AE_1$, after we made the identifications which will
ensure the associativity of the multiplication. For example, let
$A={\mathbb C}$ and let $E_i={\mathbb C}^{n_i}$ for $i=1,2$ be the
Hilbert spaces with the usual inner products and right and left
multiplications. If the isomorphism ${\mathbb C}^{n_1}\otimes
{\mathbb C}^{n_2}\rightarrow {\mathbb C}^{n_2}\otimes {\mathbb
C}^{n_1}$ is given by $e_i\otimes f_j\mapsto f_j\otimes e_i$, then
${\cal O}_E\simeq{\cal O}_{n_1}\otimes{\cal O}_{n_2}$, where
${\cal O}_n$ is the Cuntz algebra, and $\{e_1,e_2,...e_{n_1}\},
\;\{f_1,f_2,...,f_{n_2}\}$ are the canonical bases in ${\mathbb
C}^{n_i},i=1,2$.

\bigskip

{\bf Example 3}. Let $A$ be a C*-algebra, and $\alpha_1, \alpha_2$
two commuting automorphisms. The Cuntz-Pimsner algebra of the
corresponding product system is isomorphic to  the crossed product
$A\rtimes_{\alpha_1,\alpha_2}{\mathbb Z}^2$. This C*-algebra could
be obtained also from an iterated Pimsner construction. If we take
$E_i=A(\alpha_i,id)$, then ${\cal O}_{E_1}\simeq
A\rtimes_{\alpha_1}{\mathbb Z}$. By extending the scalars,
$E_2\otimes{\cal O}_{E_1}$ becomes a Hilbert bimodule over ${\cal
O}_{E_1}$, and ${\cal O}_{E_2\otimes{\cal O}_{E_1}}\simeq
A\rtimes_{\alpha_1,\alpha_2}{\mathbb Z}^2$.

\bigskip

{\bf Example 4}. For the full shift $X=\{0,1\}^{{\mathbb N}^2}$,
let $\sigma_i,\; i=1,2,$ be the horizontal and vertical shifts. In
this case, the equivalence classes in $R_n$  are homeomorphic to
the Cantor set. The existence of $\sigma_i$-systems $\mu^i$ for
$i=1,2$ will alow us to construct a product system over ${\mathbb
N}^2$, and its Cuntz-Pimsner algebra plays the role of
$C^*(\Gamma)$ defined in \S 2. More precisely, take $\alpha_i$ the
endomorphisms of $A=C(X)$ induced by $\sigma_i$, and ${\cal P}_i$
the conditional expectations  \[({\cal
P}_if)(x)=\int_{\sigma_i(y)=\sigma_i(x)}f(y)d\mu^i_x(y), i=1,2,\]
which commute. The Hilbert bimodules $A(\alpha_i,{\cal P}_i)$ are
isomorphic to $C(X)$ as vector spaces, with inner products
\[<\xi, \eta>_i(x)=\int_{\sigma_i(y)=x}\overline{\xi(y)}\eta(y)d\mu^i_x(y),\]
the left actions $ (f\cdot \xi)(x)=f(x)\xi(x)$, and the right
actions given by
\[(\xi \cdot f)(x)=\xi(x)f(\sigma_i(x)),\; i=1,2, f\in C(X).\]
The fibers $E(m,n)$ are defined as above, by taking tensor
products, and using the isomorphisms given by the Lemma. Since
none of the measures is supported on a finite set, the resulting
C*-algebra ${\cal O}_E$ is isomorphic to the Toeplitz algebra
${\cal T}_E$.

\bigskip

{\tt Valentin Deaconu, Dept of Math \& Stat/084, Univ of Nevada,
Reno, NV 89557}

{\tt  vdeaconu@unr.edu}

\end{document}